\documentclass[a4paper,10pt]{article}

\usepackage{graphicx}
\usepackage{amsfonts}
\usepackage{amstext}
\usepackage{amssymb,amsmath,amsbsy}
\usepackage[T1]{fontenc}
\usepackage[utf8]{inputenc}
\usepackage{authblk}
\usepackage{caption}

\def\R{\mathbb{R}}

\def\P{\mathbb{P}}

\def\E{\mathbb{E}}

\def\WW{\mathbb{W}}
\def\0{\mathbf{0}}
\def\1{\mathbf{1}}
\def\cc{\mathbf{c}}
\def\d{\mathbf{d}}

\def\r{\mathbf{r}}

\def\u{\mathbf{u}}

\def\M{\mathbf{M}}

\def\V{\mathbf{V}}
\def\W{\mathbf{W}}
\def\A{\mathbf{A}}
\def\B{\mathbf{B}}

\def\DD{\mathbf{D}}

\def\EE{\mathbf{E}}

\def\PP{\mathbf{P}}

\def\ep{\varepsilon}

\def\diag{\mathrm{diag}}

\def\disc{\mathrm{disc}}
\def\md{\mathrm{md}}

\def\rk{\mathrm{rank}}

\def\Vol{{\mathrm{Vol}}}

\def\part{\cal P}

\newtheorem{theorem}{Theorem}
\newtheorem{proposition}{Proposition}

\newtheorem{definition}{Definition}
\newtheorem{lemma}{Lemma}

\newtheorem{fact}{Fact}

\title{Generalized quasirandom properties of expanding graph sequences}
  
\author{Marianna Bolla}

\affil{Institute of Mathematics, Budapest University of Technology and 
Economics}

\begin{document}

\maketitle

\section*{Abstract}

We consider  multiclass spectral, discrepancy, degree, and codegree 
properties  of  expanding graph sequences. As we are able to  
prove equivalences and implications
between them and the
the definition of the generalized 
quasirandomness of Lov\'asz--S\'os (2008),  
they can be regarded as generalized quasirandom
properties akin to the equivalent
quasirandom properties of the seminal Chung--Graham--Wilson paper (1989) 
in the one-class scenario.
Since these properties are valid for certain deterministic 
graph sequences, irrespective of stochastic models, the partial implications
also justify for low-dimensional embedding of large-scale graphs and
for discrepancy  minimizing spectral clustering. 
  
\noindent
\textbf{Keywords}: generalized quasirandom graphs, multiway
discrepancy, normalized modularity spectra, cluster variances, codegrees.

\section{Introduction}\label{intro}

Our motivation comes from the multivariate statistical analysis.
The basic idea of factor analysis (static or dynamic) is to make
low-dimensional embedding of high-dimensional data. For this purpose,
the so-called $k$-factor model is investigated, where $k$ is the
hidden rank of the covariance structure and it is much less than the
number $n$ of the variables. The so-called factor scores also give rise to
partition the variables into $k$ clusters.
In the fundamental models, $k$ is fixed,
whereas $n\to \infty$ together with the sample size $N$ in the static or the  
time $T$ in the dynamic case. If our data form a graph, then $n$ is the
number of vertices, $N$ is the number of edges, and the usual condition that 
$n$ and $N$ tend to infinity 
in a prescribed way ($N$ is superlinear in $n$) 
implies that our graph is dense enough. 
Given a real-life graph on $n$ vertices, that is an
instance of an expanding sequence, with the help of graph based matrices,
our purpose is to make inference on
the number $k$ of the hidden clusters of vertices 
and to find the clusters themselves.

Graph based matrices together with eigenvalues and eigenvectors 
have been intensively studied since the 1970s.
Hoffman~\cite{Hoffman2} and Fiedler \cite{Fiedler} 
used the eigenvector, corresponding to the smallest
positive Laplacian eigenvalue of a connected graph,
to find a bipartition of the vertices which approximates the 
\textit{minimum cut problem}. 
From the two-clustering point of view ($k=2$), this eigenvector becomes
important when the corresponding eigenvalue is not separated from the trivial 
zero eigenvalue, but it is separated from the other positive
 Laplacian eigenvalues. 
On the contrary, when there is a  large spectral gap between the
trivial zero and the smallest positive Laplacian eigenvalue (or equivalently,
between the trivial 1 and the second largest positive eigenvalue of
the transition probability matrix in the random walk view), there is no use of
partitioning the vertices, the whole graph forms a highly connected cluster. 
This $k=1$ case has frequently been studied  since Cheeger~\cite{Cheeger},
establishing a lot of equivalent or near equivalent advisable features of 
these graphs.
There are many results about the relation 
between this gap and different kinds of expansion constants of the graph 
(see e.g.,~\cite{Hoory,Lovasz0}), 
including random walk view of~\cite{Azran,Diaconis,Meila}.
The vertex subsets of such graphs have a
large boundary compared to their volumes characterized
by the isoperimetric number, see Mohar~\cite{Mohar1}. 
They also show quasirandom properties discussed 
in Thomason~\cite{Thomason,Thomason1}, Bollob\'as~\cite{Bollobas}, and 
Chung, Graham, Wilson~\cite{Chung1,Chung2}.
For these favorable characteristics,
they are indispensable in communication networks.

However, less attention has been paid to graphs 
with a small spectral gap,
when several cases can occur: among others, the graph can be a bipartite
expander of Alon~\cite{Alon0} or its vertices can be divided into two sparsely 
connected
clusters, but the clusters themselves can be good expanders 
(see~\cite{Lee} and~\cite{Ng}).
In case of several clusters of vertices the situation is even more 
complicated.
The pairwise relations between the clusters and the within-cluster relations of
the vertices of the same cluster show a great variety. 
Depending on the
number and sign of the  so-called \textit{structural eigenvalues}  
of the \textit{normalized modularity matrix},
to be defined in Section~\ref{pre}, we can make inferences
on the number of the underlying clusters and the type
of connection between them.
Furthermore, based on \textit{spectral and singular value decompositions},
low dimensional embedding of the vertices is
performed, and classical and modern techniques of the
multivariate statistical analysis 
are used to find the clusters. 
The notion of the \textit{multiway discrepancy}, 
introduced in~\cite{Bolla16}, also plays a  crucial role
in identifying the clusters, and it is related to the spectral properties
of the graph.

In the worst case scenario, when there are no $k$ structural
eigenvalues with a moderate $k$, the Szemer\'edi regularity lemma
guarantees the existence of a universal $k$ (independent of $n$, it only
depends on the discrepancy bound to be attained) such that the vertices can
be classified into $k$ (equitable, and a `small' exceptional) parts such that
the between-cluster discrepancies are less than the error bound.
This theorem has overwhelming  theoretical importance, and also
spectral versions, see Szegedy~\cite{Szegedy},
whereas our purpose is to
give equivalent conditions for the existence of a $k$-cluster structure
with a moderate $k$. For this purpose, in Section~\ref{gera} we consider 
$k$-class (generalized)
random and quasirandom graphs, introduced in~\cite{LovSos}, and in
Section~\ref{gera} we establish
equivalent properties (including discrepancies and spectra) of them.
Namely, in Theorems~\ref{conj} and~\ref{ekv} we state implications and
equivalences between them. The theorems intensively use notation and 
facts of already proved theorems that are discussed in Section~\ref{pre}. 
For this long preparation, we are able to precisely formulate the main
theorems only  in Section~\ref{gera}, and prove them in Section~\ref{proofs}.
Finally, in Section~\ref{concl} we draw consequences
of our statements, as for discrepancy based spectral clustering. 

\section{Preliminaries}\label{pre}

First we recall the main result of the Chung--Graham--Wilson~\cite{Chung1}
paper about quasirandom properties that apply to the one-class quasirandomness,
when the graph sequence imitates the properties of the Erd\H os--R\'enyi
random graph with edge-density $p=\frac12$.
The authors also anticipate that instead of $\frac12$, any fixed
$0<p<1$ can be considered. Hereby, we enlist only those properties that 
will be used later for our purposes, 
and together with the original formulation (with $p=\frac12$) we give 
the analogous form with a general $p$, while we use the notation 
of the original~\cite{Chung1}  paper.

Let $(G_n )$ be a sequence of graphs as $n\to \infty$. The 
vertex-set of the general term $G_n$ is $V_n$, and
$|V_n |=n$; whereas, the number of edges of $G_n$ is $e(G_n )$.  
Consider the following properties. 
\begin{itemize}
\item
$P_1 (s)$: for all graphs $M(s)$ on $s$ vertices,
$
 {N_{G_n}^* (M(s))} = (1+o(1)) {n^s}\left( \frac12 \right)^{s \choose 2} 
$,
where $N_{G_n}^* (M(s))$ denotes the number of  
labelled induced subgraphs of $G_n$,
isomorphic to $M(s)$. With a general $p$ it reads:
$$
 {N_{G_n}^* (M(s))} = (1+o(1)) {n^s}
 {p}^{e(M(s))}{(1-p)}^{{s \choose 2} -e(M(s)) } ,
$$
where $e(M(s))$ is the number of edges in $M(s)$.

\item
$P_2 (t)$: $e(G_n ) \ge (1+o(1) ) \frac{n^2}{4}$ and
 ${N_{G_n} (C_t) } \le (1+o(1)){n^t} (\frac12 )^{t}$,
where $C_t$ is the cycle with $t$ edges and $N_{G_n} (C_t) $ is the
number of its occurrences as a (not necessarily induced) subgraph of
$G_n$, i.e., the number $\hom (C_t, G_n )$ of the $C_t \to G_n$ homomorphisms. 
(Note that a relation between $N$ and $N^*$ is given in~\cite{Chung1}.)
With a general $p$ it reads:
$$
 2e(G_n ) \ge (1+o(1) ) p n^2, \quad \textrm{and} \quad
 {\hom (C_t, G_n )}\le (1+o(1)){n^t} {p}^{t} .
$$

\item
$P_3$: $e(G_n ) \ge (1+o(1) ) \frac{n^2}{4}$, $\lambda_1 =(1+o(1))\frac{n}2$,
  $\lambda_2 =o(n)$,
where $\lambda_1$ and $\lambda_2$ are the largest and the second largest
(in absolute value) eigenvalues of the adjacency matrix of $G_n$. (Because
of the Frobenius theorem, $\lambda_1$ is always positive.)
With a general $p$ it reads:
\begin{equation}\label{*}
 2e(G_n ) \ge (1+o(1) ) p n^2, \quad
 \lambda_1= (1+o(1)) {p}n , \quad \lambda_2 =o(n) .
\end{equation}

\item
$P_4$: $\forall S\subseteq V_n$, $e(S)=\frac14 |S|^2 +o(n^2 )$, 
where $e(S)$ is the number of edges in the subgraph of $G_n$ induced by $S$.
With a general $p$ it reads: 
\begin{equation}\label{**}
 \forall X\subseteq V_n : \quad e (X,X)=p |X|^2 +o(n^2 ) ,
\end{equation}
where $e(X ,X) =2 e(X)$, and the notation $e(X, Y)$ will be used later
for the number of cut-edges between vertex-subsets $X$ and $Y$
(counting the possible edges in $X\cap Y$ twice). 


\item
$P_7$: $\sum\limits_{u,v \in V_n } |N_2 (u,v) -\frac{n}4 | =o(n^3)$.
With a general $p$ it reads: 
$$
 \sum_{u,v \in V_n } |N_2 (u,v) -p^2 n | =o(n^3) ,
$$
where $N_2 (u,v)$ is the number of common neighbors of $u,v$
in $G_n$. 
\end{itemize}
The main theorem of~\cite{Chung1} states that for $s\ge 4$ and $t\ge 4$ even,
$$
 P_2 (4) \Rightarrow P_2 (t) \Rightarrow P_1 (s) \Rightarrow
 P_3 \Rightarrow P_4 \Rightarrow \dots \Rightarrow P_7 \Rightarrow P_2 (4).
$$
The authors of~\cite{Chung1} 
called a graph \textit{quasirandom} if it satisfies any,
and therefore all, of the above (and some other, here not used)  properties.

In the $k$-class scenario, generalized quasirandom graphs imitate the
properties of the generalized random ones. 
The definition of a generalized random graph sequence is as follows
(see, e.g.,~\cite{Abbe,LovSos,SimonovitsS}) .

\begin{definition}\label{genrand}
We are given a model graph $H$ on $k$
vertices with vertex-weights $r_1 ,\dots ,r_k$ ($r_i >0$, $\sum_{i=1}^k r_i =1$)
and edge-weights $p_{ij} =
p_{ji}$, $1\le i\le j \le k$ (entries of the $k\times k$
symmetric probability matrix $\PP$, where 
$0\le p_{ij} \le 1$, $1\le i\le j \le k$, and the diagonal
entries correspond to loops at every vertex). $G_n $ is the general term
of a generalized random graph sequence on the model graph $H$ if 
\begin{itemize}
\item it has $n$ vertices;
\item to each vertex $v$ a cluster membership 
$c_v \in \{ 1, \dots ,k \}$ is assigned
according to the probability distribution $r_1 ,\dots ,r_k$;
\item
given the memberships, each pair $v\ne u$ of vertices is connected with 
probability  $p_{c_v c_u}$;  
\item further, all these decisions are  made independently.
\end{itemize}
\end{definition}
Note that in Definition~\ref{genrand}, when
$(U_1 ,\dots ,U_k)$ denotes the membership based clustering of
the vertices (they also depend on $n$, however we will not denote this
dependence, unless necessary), the following
\textbf{strong balancing condition} on the growth of the cluster sizes
can be established: with the notation
$n_i =|U_i|$, $i=1,\dots ,k$  $(\sum_{i=1}^k n_i =n)$, if
$n\to\infty$, then $\frac{n_i}{n} \to r_i$ $(i=1,\dots ,k )$.

Lov\'asz and S\'os~\cite{LovSos} gave the following definition 
of a generalized quasirandom graph sequence.
\begin{definition}\label{qgenrand}
Given a model graph $H$ on $k$
vertices with vertex-weights $r_1 ,\dots ,r_k$ and edge-weights $p_{ij} =
p_{ji}$, $1\le i\le j \le k$ (entries of $\PP$), $(G_n )$ is $H$-quasirandom
if $G_n \to W_H$ as $n\to\infty$ (in terms of the convergence of 
homomorphism densities, and $W_H$ is the step-function graphon assigned to $H$).
\end{definition}
The above graph-converge means that
for every fixed simple graph $F$, the number of copies of $F$ in $G_n $
is asymptotically the same as the number of copies of $F$ in a generalized
random graph on $n$ vertices and the same model graph $H$.
More precisely, by~\cite{Borgs},
we say that $G_n \to W_H$ if for any simple graph $F$,
$$
\frac{\hom (F, G_n )}{|V(G_n)|^{|V(F)|} } \to \hom (F, H ) 
 =\sum_{\psi : V(F) \to V(H)} \prod_{i\in V(F)} r_{\psi (i)}
 \prod_{ ij\in E(F)} p_{\psi (i) \psi (j) } .
$$
If $|V(F)| =s$, then we also have an integral formula for the
homomorphism density of the simple graph $F$ in the graphon $W_H$:
$$
\hom (F, H ) =\hom (F, W_H ) =
 \int_{[0,1]^s} \prod_{\{ i,j\} \in E(F)} W_H (x_i ,x_j ) \, dx_1 \dots dx_s .
$$
Note that the $G_n \to W_H$ convergence is also equivalent to the following:
the cut-distance between the graphons $W_{G_n}$ and $W_H$ tends to 0 as
$n\to \infty$.

As Definition~\ref{qgenrand} lays the foundation of the generalized 
quasirandomness, we want to
establish equivalent properties resembling those of the Chung-Graham-Wilson
paper. Actually, Definition~\ref{qgenrand} complies with the above $P_1 (s)$
property, since isomorphism  densities can be related to homomorphism ones.

We also introduce some
notation and former theorems (called facts) that are needed in 
our main theorems.
Let  $G=(V, \A)$ be an  undirected, \textit{edge-weighted graph}  on the
$n$-element vertex-set $V$ with the $n\times n$ symmetric weighted adjacency
matrix $\A$; the  entries satisfy $a_{ij}=a_{ji}\ge 0$, $a_{ii}=0$ and they
are similarities between the vertex-pairs. 
If we have a \textit{simple graph} $G$, then
$\A$ is the usual 0-1 adjacency matrix.

The \textit{modularity matrix}
of $G$ is defined as $\M =\A -\d \d^T $ (see~\cite{Newman}),
where the entries of $\d$ are the \textit{generalized vertex-degrees}
$d_i =\sum_{j=1}^n a_{ij}$ $(i=1,\dots ,n )$, and $\A$ is normalized in such
a way that $\sum_{i=1}^n \sum_{j=1}^n a_{ij}=1$; this assumption 
does not hurt the generality, since neither the forthcoming
\textit{normalized modularity matrix} nor the multiway discrepancies to
be defined are affected by the scaling of the entries of $\A$.
The normalized modularity matrix of $G$ (see~\cite{Bolla11}) is
$$
  \M_{\DD} =\DD^{-1/2} \M \DD^{-1/2} = \DD^{-1/2}\A\DD^{-1/2}-\sqrt{\d }\sqrt{\d }^T
  =\A_{\DD}  -\sqrt{\d }\sqrt{\d }^T ,
$$
where $\DD =\diag (d_1 ,\dots ,d_n )$ is the diagonal \textit{degree-matrix},
$\A_{\DD }$ is the normalized adjacency matrix,
and  the vector $\sqrt{\d }=(\sqrt{d_1},\dots , \sqrt{d_n})^T$ has
 unit norm.

We will assume that $G$ is \textit{connected},
i.e., $\A$ is \textit{irreducible},
in which case, the generalized vertex-degrees are all positive.
For the relation between the normalized modularity and Laplacian matrices,
see~\cite{Bolla13}.

Let $1<k<n$ be a fixed integer.
In the modularity based spectral clustering,
we look for the proper $k$-partition $U_1,\dots ,U_k $ of the vertices
such that the within- and between cluster discrepancies are minimized.
We refer to $U_i$'s as \textit{clusters}.
To motivate the introduction of the exact discrepancy measure observe that
the $ij$ entry of $\M$ is $a_{ij} - d_i d_j$, which is the difference between
the actual connection of the vertices $i,j$ and the connection that is
expected under independent attachment of them with probabilities $d_i$ and
$d_j$, respectively. Consequently, the difference between the actual
and the expected connectedness of the subsets $X,Y \subset V$ is
$$
 \sum_{i\in X} \sum_{j\in Y} (a_{ij} - d_i d_j ) =a(X,Y) -\Vol (X) \Vol (Y),
$$
where $a (X, Y) =\sum_{i\in X} \sum _{j\in Y} a_{ij}$ is the
\textit{weighted cut} between $X$ and $Y$, and
$\Vol (X) = \sum_{i\in X} d_{i}$ is the \textit{volume} of the vertex-subset $X$.
When $\A$ is the 0-1 adjacency matrix, $a(X,Y)=e(X,Y)$ is the number of 
cut-edges between $X$ and $Y$. 
Further, let $\rho (X,Y) :=\frac{a(X,Y)}{ \Vol (X) \Vol (Y)}$ be the
\textit{volume-density} between $X$ and $Y$. With these, the following
definition is given.

\begin{definition}[Definition 6 of \cite{Bolla16}]\label{diszkrepancia}
The multiway discrepancy  of the undirected, edge-weighted graph $G=(V,\A )$
in the proper $k$-partition (clustering) $U_1 ,\dots ,U_k$ of its vertices is
$$
 \md (G ; U_1 ,\dots ,U_k ) =
 \max_{1\le i \le j\le k} \max_{X\subset U_i ,\, Y\subset U_j }
 \md (X,Y;U_i , U_j ) ,
$$
where
$$
\begin{aligned}
\md (X,Y;U_i , U_j ) &=
 \frac{|a (X, Y)-\rho (U_i,U_j ) \Vol (X)\Vol (Y)|}{\sqrt{\Vol(X)\Vol(Y)}} \\
 &=| \rho (X,Y) -\rho (U_i,U_j ) | \sqrt{\Vol(X)\Vol(Y)} .
\end{aligned}
$$
The minimum $k$-way discrepancy  of  $G$ is
$$
 \md_k (G ) = \min_{(U_1 ,\dots ,U_k ) \in {\part}_k }
 \md (G ; U_1 ,\dots ,U_k ) ,
$$
where ${\part}_k$ denotes the set of proper $k$-partitions of $V$.
\end{definition}
Observe that $\md (X,Y;U_i , U_j )$ is unaffected under scaling the 
edge-weights.
Also note that $\md (G;  U_1 ,\dots ,U_k )$ is the smallest
$\alpha$ such that for every $U_i ,U_j$ pair and for every
$X\subset U_i$, $Y\subset U_j$,
$$
 |a (X, Y)-\rho (U_i,U_j ) \Vol (X)\Vol (Y)| \le \alpha \sqrt{\Vol(X)\Vol(Y)}
$$
holds.
It resembles the notion of volume-regular cluster pairs of~\cite{Alon10}
or the $\epsilon$-regular pairs in the Szemer\'edi
regularity  lemma~\cite{Szemeredi}, albeit with given number of
vertex-clusters, which are usually not equitable;
further, with volumes, instead of cardinalities.

The forthcoming Facts~1 and~2 justify  for
the following spectral relaxation of the minimum $k$-way
discrepancy problem.
Let the eigenvalues of $\M_{\DD}$, enumerated in decreasing
absolute values, be $1\ge |\mu_1| \ge |\mu_2 | \ge\dots \ge |\mu_n |=0$.
Assume that $|\mu_{k-1}| > |\mu_k|$, and denote by $\u_1 ,\dots ,\u_{k-1}$
the unit-norm, pairwise orthogonal eigenvectors, corresponding to
$\mu_1 ,\dots ,\mu_{k-1}$.
Let $\r_1 ,\dots ,\r_n \in \R^{k-1}$ be the row vectors of the $n\times (k-1)$
matrix of column vectors $\DD^{-1/2} \u_1 ,\dots , \DD^{-1/2} \u_{k-1}$;
they are called $(k-1)$-dimensional representatives of the vertices.

The \textit{weighted k-variance} of these representatives is defined as
\begin{equation}\label{wkszoras}
 {\tilde S}_k^2  =\min_{(U_1 ,\dots ,U_k ) \in {\part}_k }
\sum_{i=1}^k \sum_{v\in U_i } d_v \| \r_v -{ \cc }_i \|^2 ,
\end{equation}
where ${\cc }_i =\frac1{\Vol (U_i ) } \sum_{v\in U_i } d_v \r_v $ is the
weighted center of the cluster $U_i$.
It is the \textit{weighted k-means algorithm} that provides this minimum.
We will also need the plain $k$-variance of the
representatives $\r_1 ,\dots ,\r_n \in \R^k$ that are row-vectors of the
matrix, the columns of which are the unit-norm, pairwise orthogonal
eigenvectors corresponding to the $k$ largest (in absolute value) eigenvalues
of $\A$.  This \textit{k-variance} is
\begin{equation}\label{kszoras}
 {S}_k^2  =\min_{(U_1 ,\dots ,U_k )}
\sum_{i=1}^k \sum_{v\in U_i } \| \r_v -{ \cc }_i \|^2 ,
\end{equation}
where ${\cc }_i =\frac1{ |U_i | } \sum_{v\in U_i } \r_v $ is the
center of the cluster $U_i$.
It is the usual \textit{k-means algorithm} that finds this minimum.
By an easy analysis of variance argument it follows that the 
optimum ${S}_k$ is just the minimum  distance
between the subspace spanned by the eigenvectors corresponding to the
$k$ largest (in absolute value) eigenvalues of $\A$   and  the one
of the step-vectors over the $k$-partitions of $V$.
Similar holds for ${\tilde S}_k$ with the transformed eigenvectors.

Note that finding a global minimizer for the $k$-means problem is NP-complete.
However, there are efficient polynomial time algorithms for finding an
approximate solution whose value is within a constant fraction of the
optimal value,  under certain conditions. In particular,
Theorem 4.6 of~\cite{Ostrovsky} states that if the data
satisfy the $k$-clusterable criterion ($S_k^2 \le \epsilon^2 S_{k-1}^2$
with a small enough $\epsilon$), then there is a PTAS (polynomial time
approximation scheme) for the $k$-means problem.
Our conditions for the $k$-variances in Theorems~\ref{conj} and~\ref{ekv} 
do comply with this requirement.
 
The spectral relaxation means that we can approximately find discrepancy
minimizing clustering via applying the unweighted or weighted $k$-means 
algorithm
to the $k$- or $(k-1)$-dimensional vertex representatives. 
This is supported by the following facts.

\begin{fact}[Theorem 7 of \cite{Bolla16}]\label{undirect}
Let $G= (V, \A)$ be an edge-weighted, undirected graph, $\A$ is
irreducible. Then for any integer $1\le k < \rk (\A)$, 
$$
 |\mu_k |  \le 9\md_{k } (G )  (k+2 -9k\ln \md_{k } (G )) 
$$
holds, provided $0< \md_{k } (G ) <1$,
where $\mu_k$ is the $k$-th largest eigenvalue (in absolute value) of
the normalized modularity matrix $\M_{\DD}$ of $G$. 
\end{fact}
The above Fact~\ref{undirect} is a certain converse of the multiclass
expander mixing lemma.  
In the forward direction, in the $k=1$ case, 
the following result is considered as the
extension of the expander mixing lemma to irregular graphs:
$$
 \md_1 (G) \le \| \M_{\DD} \| =|\mu_1 |,
$$
where $\| \M_{\DD} \|$ is the spectral norm of the normalized modularity matrix
of $G$. This result was proved in~\cite{Chung2}, where
$\md_1 (G)$ is denoted by $\disc (G)$.
When $k\ge 1$, we are able to prove the following generalization of the
expander mixing lemma.

\begin{fact}[Theorem 3 of~\cite{Bolla17}]\label{thkuj}
Let $G_n $ be the general term of a connected simple graph sequence, $G_n$ has
$n$ vertices. (We do not denote the dependence of the vertex-set $V$ and
adjacency matrix $\A$ of $G_n$ on $n$). 
Assume that there are constants $0<c<C<1$ 
such that except $o (n)$ vertices, the degrees satisfy 
$cn \le d_v \le Cn$, $v=1,\dots ,n$.   
Let the eigenvalues  of the normalized modularity matrix $\M_{\DD}$ of $G_n$,
enumerated in decreasing absolute values, be 
$$
 |\mu_1 | \ge \dots \ge |\mu_{k-1} | >\ep \ge |\mu_k | \ge
 \dots \ge |\mu_n |=0 . 
$$
The partition $(U_1,\dots ,U_k )$ of $V$ is defined so that it minimizes
the weighted k-variance $s^2 ={\tilde S}_k^2$ 
of the optimal $(k-1)$-dimensional vertex representatives of $G_n$.
Assume that $(U_1, \dots ,U_k )$ satisfies the strong balancing condition.
Then 
$$
 \md (G_n ; U_1 ,\dots ,U_k ) \le 2 
  \left( \frac{C}{c} +o(1) \right) (\sqrt{2k} s +\ep ) .
$$
\end{fact}
Fact~\ref{thkuj} implies that 
$\md_k (G ) \le  \md (G ; U_1 ,\dots ,U_k )= 
{\cal O} (\sqrt{2k} {\tilde S}_k +|\mu_k | )$. For the $k=1$ case, 
$\ep =|\mu_1 |$,
${\tilde S}_1 =0$ (based on the coordinates of the 
$\DD^{-1/2} \u_1 =\DD^{-1/2} \sqrt{\d } =\1$ vector),
and so, we get back the original expander mixing lemma up to a constant.
In the $k=2$ bipartite, biregular case we get the statement of~\cite{Evra},
see~\cite{Bolla16} for further explanation.
Consequently, 
a `small' $|\mu_k |$ and ${\tilde S}_k$ is an indication of $k$
clusters with `small' within- and between-cluster discrepancies;
further, a partition, close to the optimum, can be obtained by spectral tools.
Also observe that the degree-condition of Fact~\ref{thkuj} means that 
our graph is dense  enough.

\section{Generalized random and quasirandom properties with 
equivalence theorems}\label{gera}

Now some properties  of a generalized random graph sequence are stated
so that to get ideas how to formulate generalized quasirandom properties.

\begin{proposition}\label{ekvgen}
Let $(G_n )$ be a generalized random graph sequence  on the model graph $H$;
$G_n$ has $n$ vertices with
vertex-classes  $U_1 ,\dots ,U_k $ of sizes $n_1 ,\dots n_k$
(they also depend on $n$).
Let $H$, and so $k$ be kept fixed,
i.e., the $k\times k$ 
probability matrix $\PP$ of rank $k$ 
and the `blow-up' ratios $r_1 ,\dots ,r_k$ are 
fixed, while 
$n\to\infty$  under the strong balancing condition.
Then the following properties hold \textit{almost surely}
for the convergence of the graph sequence  $(G_n )$, for the
adjacency matrix $\A_n =(a_{ij}^{(n)})$, the normalized modularity matrix 
$\M_{\DD,n}$, the multiway discrepancies, and the within- and 
between-cluster codegrees  of $G_n$.
\begin{itemize}
\item[{{0.}}]
$G_n \to W_H$ as $n\to\infty$.

\item[{{1.}}]
$\A_n$ has exactly $k$ so-called \textit{structural} eigenvalues
that are $\Theta (n)$, 
while the remaining eigenvalues are $O(\sqrt{n})$ (in absolute value).
Further, the $k$-variance $S_{k,n}^2$ (see (\ref{kszoras})) of the
$k$-dimensional vertex representatives, based on the eigenvectors corresponding
to the structural eigenvalues of $\A_n$, is $O(\frac1{n})$.

\item[{{2.}}]
There exists a positive constant $0<\delta <1$ independent of $n$ (it only
depends on $k$) such that $\M_{\DD,n}$ has
exactly $k-1$ \textit{structural} eigenvalues of 
absolute value greater than $\delta$,
while all the other eigenvalues are  $O (n^{-\tau})$ 
for every $0<\tau < \frac12$.
Further, the  weighted $k$-variance ${\tilde S}_{k,n}^2$ (see (\ref{wkszoras}))
of the $(k-1)$-dimensional vertex representatives, based on the transformed
eigenvectors corresponding
to the structural eigenvalues of $\M_{\DD,n}$,
is $O(n^{-2\tau})$,  for every $0<\tau < \frac12$.

\item[{{3.}}]
There is a constant $0<\theta <1$ (independent of $n$)
such that $\md_1 (G_n )>\theta$, \dots ,
$\md_{k-1} (G_n ) >\theta$, and
the $k$-way discrepancy $\md (G_n ; U_1 ,\dots ,U_k  )$ 
is $O(n^{-\tau})$,  for every $0<\tau < \frac12$. 
 
\item[{{4.}}]
For every $1\le i\le j\le k$ and $u \in U_i$: 
$$
 N_1 (u; U_j ):= \sum_{v\in U_j} a_{uv}^{(n)} =p_{ij} n_j +o(n)  ,
$$
where $N_1 (u; U_j )$ denotes the number of neighbors of $u$ in $U_j$.

For every $1\le i\le j\le k$ and $u,v \in U_i$, $u\ne v$: 
$$
 N_2 (u,v; U_j):=\sum_{t \in U_j} a_{ut }^{(n)}  a_{vt }^{(n)} =p^2_{ij} n_j +o(n) ,
$$
where  $N_2 (u,v; U_j)$ denotes the number of common neighbors of $u,v$ in
$U_j$. 
\end{itemize}
\end{proposition}

By Property 0, a generalized random graph sequence is also
generalized quasirandom, though it converges more quickly. 
Property~0 was proved in~\cite{Bolla12}, while
Property~1 in~\cite{Bolla5}, and Property~2 in~\cite{Bolla10}
in a more general framework of the SVD of
rectangular arrays of nonnegative entries,
they are summarized in~\cite{Bolla13}. 
The basic idea of the proofs is that the adjacency matrix $\A_n$ of  $G_n$
can be decomposed as a deterministic
block-matrix $\B_n$ (blown-up of $\PP$ with sizes $n_1 ,\dots ,n_k$) 
plus a so-called Wigner-noise $\W_n$, see~\cite{Bolla5}. 
The Wigner-noise has spectral norm
$O(\sqrt{n})$ almost surely, 
while the blown-up matrix has as many non-zero eigenvalues of order
$n$ as the rank of $\PP$ (in our case, $k$) with eigenvectors  
that are stepwise constant over $U_1 ,\dots ,U_k$.
 
Note that to prove Properties
1 and 2, even the following
\textbf{weak balancing condition} suffices:
$n\to\infty$ in such a way that $\frac{n_i}{n} \ge c$ $(i=1,\dots ,k )$
with some constant $0<c\le \frac1{k}$.
Property~3 is the consequence of Property~2, by the
the back and forth statements 
between discrepancy and normalized modularity spectra, see 
Facts~\ref{undirect} and~\ref{thkuj}.  

Property~4 is easy to prove by large deviations.
The subgraph of $G_n$,
induced by $U_i$ and denoted by $G_{ii,n}$, is the general term of
an Erd\H os--R\'enyi type
random graph sequence with edge probability $p_{ii}$, for every $i=1,\dots ,k$.
The bipartite subgraph of $G_n$,
induced by the $U_i ,U_j$ pair  and denoted by $G_{ij,n}$, is the general term
of a bipartite random graph sequence with edge probability $p_{ij}$, for
every $i,j=1,\dots ,k$; $i\ne j$ pair. Therefore, the subgraphs are
almost surely regular, while the bipartite subgraphs are almost surely 
biregular. 
Consequently, the vertex-degrees are of order $\Theta (n)$, almost surely.
Further, the codegrees are as expected: every two vertices in $U_i$ have 
approximately the same number of common neighbors in $U_j$, 
for every $i,j=1,\dots ,k$.

Note that the above generalized random graph in another context is discussed as
the stochastic block model or planted partition model, 
see, e.g.,~\cite{Cojab,Holland,McSherry}, 
though these papers work with a fixed $n$ and
do not consider any condition for the growth of the cluster sizes. 
Indeed, Definition~\ref{genrand} provides us with a random graph model
without the hidden (planted) clusters revealed. For this purpose, there are
algorithms available, e.g., in~\cite{Cojab,McSherry};
however, one wonders whether a large and dense
enough real-life graph can be `close' to a one coming from this model. 
In the sequel, we will define precisely some properties that 
are weaker than those of Proposition~\ref{ekvgen}, but can
characterize a class of graphs, given $k$. 
These will be called generalized quasirandom properties. 
Note that some other papers, e.g., \cite{Abbe,Bollobas3,Mossel} 
scale the probability
matrix with $n$, and prove the consistency of the clustering under these
conditions.

Properties, reminiscent of those of the generalized random graphs,
are now  formulated 
for expanding deterministic graph sequences, and we show that there are
many equivalences and implications between them, 
irrespective of stochastic models.
Theorem~\ref{conj}
states mainly implications, whereas Theorem~\ref{ekv} states equivalences.

\begin{theorem}\label{conj}
Let $G_n $ be the general term of a sequence of simple graphs  with
vertex-set $V_n$, adjacency matrix $\A_n =(a_{ij}^{(n)})$,
and normalized modularity matrix $\M_{\DD,n}$.
Let $k$ be a fixed positive integer, whereas $|V_n |=n \to \infty$.
Consider the following properties:
\begin{itemize}
\item[{{P0.}}]
There exists a vertex- and edge-weighted  graph $H$ on $k$ vertices
with vertex-weights $r_1 ,\dots ,r_k$ and edge-weights $p_{ij} =
p_{ji} \in [0,1]$, $1\le i\le j \le k$, where the $k\times k$
symmetric probability  matrix $\PP =(p_{ij})$ has rank $k$,
such that $G_n \to W_H$ as $n\to \infty$. 

\item[{{PI.}}]
$\A_n$ has $k$ \textit{structural} eigenvalues $\lambda_{1,n}, \dots ,
\lambda_{k,n}$ such that the normalized eigenvalues converge:
$\frac1{n} \lambda_{i,n}  \to q_i$ as $n\to\infty$ $(i=1,\dots ,k )$
with some non-zero reals $q_1 ,\dots ,q_k$, and the
remaining eigenvalues are $o (n)$.

The $k$-variance $S_{k,n}^2$ of the
$k$-dimensional vertex representatives, based on the eigenvectors corresponding
to the structural eigenvalues of $\A_n$,  is $o(1)$.
The $k$-partition $(U_1 ,\dots ,U_k )$ minimizing this
$k$-variance satisfies the strong balancing condition.

\item[{{PII.}}]
$G_n$ has no dominant vertices: there are constants 
$0<c<C<1$ such that the vertex-degrees are between $cn$ and $Cn$,
except of possibly $o(n)$ vertices; further, there exists
a constant $0<\delta <1$ (independent of $n$) such that
$\M_{\DD,n}$ has  $k-1$ \textit{structural} eigenvalues
that are greater than $\delta$ (in absolute
value), while the remaining eigenvalues are $o(1)$.
 

The weighted $k$-variance ${\tilde S}_{k,n}^2$ of the
$(k-1)$-dimensional vertex representatives, based on the 
transformed eigenvectors corresponding
to the structural eigenvalues of $\M_{\DD,n}$,  is $o(1)$.
The $k$-partition $(U_1 ,\dots ,U_k )$ minimizing this
$k$-variance satisfies the strong balancing condition.

\item[{{PIII.}}]
There are vertex-classes  $U_1 ,\dots ,U_k $ obeying the strong
balancing condition, and
there is a constant $0<\theta <1$ (independent of $n$)
such that $\md_1 (G_n )>\theta ,\dots , \md_{k-1} (G_n  ) >\theta$, and
$\md (G_n ; U_1 ,\dots ,U_k )= o(1)$.

\item[{{PIV.}}]
There are vertex-classes  $U_1 ,\dots ,U_k $ of sizes $n_1 ,\dots ,n_k$
obeying the strong balancing condition, 
and there is a $k\times k$
symmetric  probability matrix $\PP =(p_{ij} )$ of rank $k$ such that, 
with them, the following holds:
\begin{equation}\label{cod}
 \sum_{u,v\in U_i } | N_2 (u,v; U_j ) -p_{ij}^2 n_j | = o(p_{ij}^2 n_i^2 n_j)
 = o(n^3 ), \quad \forall i,j=1,\dots ,k.
\end{equation}
\end{itemize}
Then P0 is equivalent to PIV, and they imply PI and PII; further, PII implies
PIII.
We also  consider the following strengthening of property PI:
\begin{itemize}
\item[{{PI+.}}]
$\A_n$ has $k$ \textit{structural} eigenvalues $\lambda_{1,n}, \dots ,
\lambda_{k,n}$ such that the normalized eigenvalues converge:
$\frac1{n} \lambda_{i,n}  \to q_i$ as $n\to\infty$ $(i=1,\dots ,k )$
with some non-zero reals $q_1 ,\dots ,q_k$, and the
remaining eigenvalues are $o (\sqrt{n})$.

The $k$-variance $S_{k,n}^2$ of the
$k$-dimensional vertex representatives, based on the eigenvectors corresponding
to the structural eigenvalues of $\A_n$,  is $o(\frac1{n})$.
The $k$-partition
$P_{k,n} =(U_{1n}, \dots ,U_{kn} )$ of the vertices of $G_n$ 
minimizing this $k$-variance not only satisfies the strong balancing condition, 
but we also assume that there is a $k\times k$ symmetric 
probability matrix $\PP =(p_{ij})$ of rank $k$ such that
\begin{equation}\label{edconv}
 d (U_{in},U_{jn} ) :=\frac{e(U_{in} ,U_{jn} )}{ |U_{in}| |U_{jn} | } \to p_{ij} 
 \quad (1\le i\le j\le k),  \quad n\to\infty . 
\end{equation}
(I.e., the within- and between-cluster edge densities converge to the entries
of $\PP$.)
\end{itemize}
Then PI+ implies P0.
\end{theorem}

Note that properties P0--PIV are in many ways weaker than those of the
generalized random graphs.
For example, the generalized quasirandom graphs may converge much more slowly 
than the
generalized random ones that can be decomposed as a deterministic
block-matrix plus a  Wigner-noise. The Wigner-noise has spectral norm
$O(\sqrt{n})$, therefore converges very quickly. In case of a generalized
quasirandom graph, the noise has spectral norm $o( n)$ only.
Also, in the generalized random case, 
the codegree condition holds for any pair of the vertices
(of the same cluster toward those of a given own or other cluster),
whereas in the general quasirandom case, it only holds on average,
see~(\ref{cod}). However, PI+ is even stronger than Property~1 of
Theorem~\ref{ekvgen}: it requires $o (\sqrt{n} )$ magnitude of the
non-structural eigenvalues, whereas those of a generalized random graph are
of order $O (\sqrt{n})$. However, except of $o(n)$ ones, the latter are also of 
order  $o (\sqrt{n})$, in view of
concentration results on the eigenvalues, see, e.g.,~\cite{Alon1}.
In summary, the above P0 is stronger than PI, but weaker than PI+. In the
$k=1$ case, property $P_3$ of the Chung--Graham--Wilson~\cite{Chung1}
paper is equivalent to $P_1 (s)$ for every $s$ (our P0).
But there, above the separation in the spectrum, a lower bound for the
edge-density should be satisfied. In the $k>1$ scenario,
the inter- and intra-cluster edge-densities can be bounded only
at the expense of each other, so we have to require the convergence of
these edge-densities too.

\begin{theorem}\label{ekv}
Let us define PII+ and PIII+ as PII and PIII of Theorem~\ref{conj} 
together with the following
additional assumptions, respectively.
The $k$-partition $(U_1 ,\dots ,U_k )$ emerging in PII and PIII 
of Theorem~\ref{conj} not only
satisfies the strong balancing condition, but there is a $k\times k$ symmetric 
probability matrix $\PP =(p_{ij})$ of rank $k$ such that
\begin{equation}\label{densi}
 d (U_{i},U_{j} ) = p_{ij} +o(1) 
 \quad (1\le i\le j\le k),  \quad n\to\infty  
\end{equation}
(the same as~(\ref{edconv})),
and for every $1\le i \le j \le k$ and $u\in U_i$,
\begin{equation}\label{degri}
 N_1 (u; U_j ) =(1+o(1)) p_{ij} n_j 
\end{equation}
holds. 
Then P0 $\Longrightarrow$ PII+ $\Longrightarrow$ PIII+ $\Longrightarrow$ PIV 
$\Longrightarrow$ P0, so they are all equivalent. 
\end{theorem}

\section{Proofs of Theorems~\ref{conj} and~\ref{ekv}}\label{proofs}

\noindent
\textbf{Proof of Theorem~\ref{conj}}.

\noindent
\textbf{Proof of P0 $\Longrightarrow$ PIV}: 
we use  the results of~\cite{Chung1,LovSos}.
By~\cite{LovSos}, 
the vertex set of the generalized
quasirandom graph $G_n$ (defined by P0) can be partitioned into classes 
$U_1 ,\dots ,U_k$
in such a way that $\frac{|U_i |}{n} \to r_i$ $(i=1,\dots ,k)$, 
that gives the strong balancing;
the subgraph $G_{ii,n}$  is the general term of a quasirandom
graph sequence with edge-density tending to
$p_{ii}$ $(i=1,\dots,k )$, whereas $G_{ij,n}$
is the general term of a bipartite quasirandom 
graph sequence with edge-density tending to $p_{ij}$ $(i\ne j)$
as $n\to\infty$.  
Therefore, for the subgraphs the equivalent statements of~\cite{Chung1}
of the usual (one-class) quasirandomness are applicable, and similar
considerations can be made for the bipartite subgraphs as well,
see~\cite{Bollobas,Kohayakawa,Thomason1}.
We also need two simple lemmas.

\begin{lemma}\label{l1}
If $(G_{ii,n})$ is quasirandom, then
$$
\sum_{u,v\in U_i} N_2 (u,v; U_i ) \ge (1+o(1)) p_{ii}^2 n_i^3 , 
  \quad i=1,\dots ,k. 
$$
\end{lemma}
\noindent
\textbf{Proof of Lemma~\ref{l1}:} We drop the index $n$ of the adjacency 
entries.
$$
\begin{aligned}
&\sum_{u,v\in U_i} N_2 (u,v; U_i ) =\sum_{u,v\in U_i} \sum_{t\in U_i }
a_{ut} a_{vt}= \sum_{t\in U_i } \sum_{u\in U_i } a_{ut} \sum_{v\in U_i } a_{vt}
 = \sum_{t\in U_i } [N_1 (t; U_i )]^2   \\
&\ge \frac1{n_i} 
 [\sum_{t\in U_i } N_1 (t; U_i )]^2 = \frac1{n_i} [2e(U_i )]^2 \ge \frac1{n_i}
 [(1+o(1)) p_{ii} n_i^2 ]^2 = (1+o(1)) p_{ii}^2 n_i^3  ,
\end{aligned}
$$
where  $e(U_i)$ is the number of edges within the  subgraph $G_{ii,n}$
of $G_n$, induced by $U_i$. 
In the first inequality we used the Cauchy--Schwarz, and in the second
one, the first part of the equivalent quasirandom 
property $P_3$ of~\cite{Chung1}, see~(\ref{*}). $\square$

\begin{lemma}\label{l2}
If $(G_{ij,n})$ is bipartite quasirandom, then
$$
\sum_{u,v\in U_i} N_2 (u,v; U_j  ) \ge (1+o(1)) p_{ij}^2 n_i^2 n_j , 
  \quad i\ne j. 
$$
\end{lemma}
\noindent
\textbf{Proof of Lemma~\ref{l2}:} We drop the index $n$ of the adjacency 
entries.
$$
\begin{aligned}
&\sum_{u,v\in U_i} N_2 (u,v; U_j ) =\sum_{u,v\in U_i} \sum_{t\in U_j }
a_{ut} a_{vt}= \sum_{t\in U_j } \sum_{u\in U_i } a_{ut} \sum_{v\in U_i } a_{vt}
 = \sum_{t\in U_j } [N_1 (t; U_i )]^2   \\
&\ge \frac1{n_j} 
 [\sum_{t\in U_j} N_1 (t; U_i )]^2 = \frac1{n_j} [e(U_i ,U_j )]^2 \ge 
\frac1{n_j} [(1+o(1)) p_{ij} n_i n_j ]^2 = (1+o(1)) p_{ij}^2 n_i^2 n_j  ,
\end{aligned}
$$
where recall that 
$e(U_i ,U_j )$ is the number of cut-edges between $U_i$ and $U_j$, i.e.,
the number of edges in the bipartite subgraph $G_{ij,n}$ of $G_n$,
induced by the $U_i , U_j$ pair.
Here, in the first inequality we used the Cauchy--Schwarz, and in the second
one, the equivalent quasirandom 
property of bipartite quasirandom graphs. $\square$

In view of the lemmas, we estimate the square of the left-hand side
of~(\ref{cod})  by the Cauchy--Schwarz inequality
for every $i,j=1,\dots k$, with $n_i^2$ terms: 
$$ 
\begin{aligned}
&\left\{ \sum_{u,v\in U_i } | N_2 (u,v; U_j ) -p_{ij}^2 n_j | \right\}^2 
\le n_i^2 \sum_{u,v\in U_i } | N_2 (u,v; U_j ) -p_{ij}^2 n_j |^2  \\
&= n_i^2 \left\{ \sum_{u,v\in U_i } [ N_2 (u,v; U_j ) ]^2 -2 p_{ij}^2 n_j
 \sum_{u,v\in U_i }  N_2 (u,v; U_j ) + n_i^2 (p_{ij}^2 n_j )^2 \right\}  \\
&\le n_i^2 \left\{ (1+o(1)) p_{ij}^4 n_i^2 n_j^2 -2 (1+o(1))  p_{ij}^4 n_i^2 n_j^2
+ p_{ij}^4 n_i^2 n_j^2 \right\}  \\
&=  n_i^2 o(1) p_{ij}^4 n_i^2 n_j ^2 = o( p_{ij}^4 n_i^4 n_j^2 ) ,
\end{aligned}
$$
where to estimate $\sum_{u,v\in U_i} N_2 (u,v; U_j  ) $ we used Lemma~\ref{l1}
in the $i=j$, and Lemma~\ref{l2} in the $i\ne j$ case; further,
utilized that $\sum_{u,v\in U_i} [N_2 (u,v; U_j  )]^2 $ is
asymptotically $\hom (C_4 , G_{ij,n} )$, the number of the 
$C_4 \to G_{ij,n}$ homomorphisms, where $C_4$ is the 4-cycle graph.
Indeed, in the $i=j$ case, the equivalent property $P_2(4)$ of~\cite{Chung1}
guarantees that $\hom (C_4 , G_{ii,n} ) \le (1+o(1)) p_{ii}^4 n_i^4$, the latter
is  what is expected in an Erd\H os--R\'enyi 
type random graph with edge-density $p_{ii}$.
In the $i\ne j$ case a similar property (see~\cite{LovSos}) for
bipartite quasirandom graphs implies that the $C_4 \to G_{ij,n}$
homomorphism density is 
$$
 \frac{\hom (C_4 , G_{ij,n} )}{n_i^2 n_j^2} = (1+o(1)) p_{ij}^4 .
$$
Here $\hom (C_4 , G_{ij,n} ) = \sum_{u,v\in U_i} [N_2 (u,v; U_j  )]^2 $
asymptotically,
as by~\cite{LovSos} only 4-cycles in the above bipartition have to be 
considered; these 4-cycles have two vertices from $U_i$ and two from 
$U_j$, and any two of the
common neighbors of $u,v\in U_i$ in $U_j$ are possible candidates
to close a (labelled) 4-cycle with them.$\square$

\vskip0.5cm
\noindent
\textbf{Proof of PIV $\Longrightarrow$ P0}: 
If the average codegree condition~(\ref{cod}) holds
for the subgraphs $(i=j)$, then by the $P_7 \Rightarrow P_1(s)$ $(s=1,2,\dots )$
implication of the equivalent quasirandom properties of~\cite{Chung1},
the subgraphs $G_{ii,n}$ are quasirandom (in terms of the 
isomorphisms, and so of the homomorphism
densities). Likewise, in the $i\ne j$ case, the bipartite subgraphs
$G_{ij,n}$ are bipartite quasirandom. 
Therefore, $G_n$ is built of quasirandom  and bipartite quasirandom blocks, 
so under the strong
balancing condition, they together form a generalized
quasirandom graph sequence on $k$ classes and model graph $H$, the
vertex-weights of which are $r_1 ,\dots ,r_k$ of the strong balancing
condition, and the
edge-weights are entries in the
condition~(\ref{cod}) of the probability matrix $\PP$ in PIV. $\square$

\vskip0.5cm
\noindent
\textbf{Proof of P0 $\Longrightarrow$ PII:} 
We will use the following results proved in former papers.

\begin{fact}[Theorem 8 of \cite{Bolla14}]\label{sajkonv}
Let $G_n =(V_n , \A_n )$ be the general term of a convergent sequence of 
connected edge-weighted graphs whose edge-weights are in [0,1] and 
the vertex-weights are the generalized degrees. Assume that
there are no dominant vertices. Let $W$ denote the limit graphon of the
sequence $(G_n )$, and let
$$
 |\mu_{n,1}  |\ge  |\mu_{n,2}| \ge \dots \ge |\mu_{n,n}| =0
$$
be the normalized modularity spectrum of $G_n$ (the eigenvalues are indexed
by their  decreasing absolute values).
Further, let $\mu_i (P_{\WW })$ be the $i$-th largest (in absolute value) 
eigenvalue of the  integral operator 
$P_{\WW } : L^2 (\xi' ) \to L^2 (\xi )$ 
taking conditional expectation with respect to
the joint measure $\WW$ embodied by the normalized limit graphon $W$,
and $\xi ,\xi'$ are identically distributed random variables
with the marginal distribution of their
symmetric joint distribution  $\WW$, whereas $L^2 (\xi )$ denotes the
Hilbert space of the measurable functions of $\xi$ with 0 expectation and
finite variance.
Then for every $i\ge 1$, 
$$
  \mu_{n,i} \to \mu_i (P_{\WW }) \quad \textrm{as} \quad  n\to\infty .
$$
\end{fact}

\begin{fact}[Theorem 9 of \cite{Bolla14}]\label{alterkonv}
Assume that there are constants $0<\ep <\delta \le 1$ such that
the normalized modularity spectrum  of the above $G_n$  satisfies 
$$
 |\mu_{n,1}  |\ge \dots \ge |\mu_{n,k-1}| \ge \delta >\ep \ge 
|\mu_{n,k}| \ge \dots \ge |\mu_{n,n} | =0 .
$$
With the notation of Fact~\ref{sajkonv}, and assuming that there are
no dominant vertices of $G_n$, the subspace spanned by the vectors
$\DD_n^{-1/2}\u_{n,1}$, \dots ,$\DD_n^{-1/2}\u_{n,k-1}$, where
$\u_{n,1}, \dots ,\u_{n,k-1}$ are orthonormal eigenvectors
 belonging to the $k-1$ largest absolute value
eigenvalues of the normalized modularity matrix of $G_n$, also converges to
the corresponding $(k-1)$-dimensional subspace of $P_{\WW}$. 
More exactly, if $\PP_{n,k-1}$ denotes
the projection onto the subspace spanned by the vectors
$\DD_n^{-1/2}\u_{n,1}$, \dots ,$\DD_n^{-1/2}\u_{n,k-1}$,
and $\PP_{k-1}$ denotes the projection 
onto the  analogous eigen-subspace of $P_{\WW}$, then 
$\|  \PP_{n,k-1} - \PP_{k-1} \| \to 0$ as $n\to\infty $
(in spectral norm).
\end{fact}

So the proof of P0 $\Longrightarrow$ PII is as follows.
Since by P0, $(G_n)$ converges to the limit graphon $W_H$, and
the eigenvalues $\mu_i (P_{\WW_H })$ of Fact~\ref{sajkonv} consist of $k-1$
non-zero numbers, and the others are zeros, the statement for the convergence
of the modularity spectrum follows. As for the 
weighted $k$-variances, we use
Fact~\ref{alterkonv}, which implies that the subspace spanned by the
transformed eigenvectors corresponding to the structural eigenvalues of 
$\M_{\DD,n}$ converges to the subspace of step-vectors; further, the steps
are proportional to $r_i$'s, so the strong balancing also follows.
As the weighted $k$-variance depends continuously on the above subspaces,
Fact~\ref{alterkonv} implies the convergence of the 
weighted $k$-variance as well. 
Note that here $\| \PP_{n,k-1} - \PP_{k-1} \|_F \le \sqrt{k-1} 
\| \PP_{n,k-1} - \PP_{k-1} \| \to 0$ as $n\to \infty$,
where $\| \cdot \|_F$ denotes the Frobenius norm. $\square$ 

\vskip0.5cm
\noindent
\textbf{Proof of P0 $\Longrightarrow$ PI}:
We use Theorem 6.7 of~\cite{Borgs3}, where the authors  prove
that if the sequence $(W_{G_n})$ of graphons converges  to the limit 
graphon $W$,
then both ends of the spectra of the integral operators, 
induced by $W_{G_n}$'s 
as kernels, converge to the ends of
the spectrum of the integral operator induced by $W$ as 
kernel. We apply this argument to the step-function graphon corresponding to
$G_n$ (the eigenvalues of the induced integral operator are the normalized
eigenvalues of $G_n$) and for the limit graphon $W_H$ of $(G_n )$.
The same argument as in P0 $\Longrightarrow$ PII 
can be applied for the convergence of the spectral subspaces,
so by the above considerations,
the convergence of the $k$-variances is also obtained.
Since the steps of the emerging step-vectors
are proportional to $r_i$'s, the strong balancing condition also holds.
$\square$

\vskip0.5cm
\noindent
\textbf{Proof of PII $\Longrightarrow$ PIII}: We will use Facts~\ref{undirect}
and~\ref{thkuj}.

Assume that there is a constant $0<\delta <1$  such that
$\M_{\DD,n}$ has  $k-1$  eigenvalues
that are greater than $\delta$ in absolute
value, while the remaining eigenvalues are $o(1)$;
further, the squareroot of weighted $k$-variance ${\tilde S}_{k,n}^2$ is 
also $o(1)$. Using that there are no dominant vertices,
we apply Fact~\ref{thkuj}. 
According to this, 
$\md (G_n ; U_1 ,\dots ,U_k )= o (1)$. Indirectly, assume that there is no 
absolute constant $0<\theta <1$ 
such that $\md_1 (G_n ) >\theta ,\dots , \md_{k-1} (G_n  ) >\theta$. 
Then there is
an $1\le i\le k-1$ with $\md_i (G_n ) \le \ep$ for any $0<\ep <1$.
But  Fact~\ref{undirect} estimates $|\mu_{n,i}|$
with a (near zero) strictly increasing function of $\md_i (G_n )$. In view
of this, there should be an $0<\ep' <1$ so that $|\mu_{n,i}| \le \ep'$, where
$\ep'$ can be any small positive number (depending on $\ep$). 
This contradicts to the  $|\mu_{n,i}| > \delta$ assumption. $\square$

\vskip0.5cm
\noindent
\textbf{Proof of PI+ $\Longrightarrow$ P0:} 
Under the assumptions of PI+,
by Theorem 3.1.17 of~\cite{Bolla13} we are able to find a blown-up matrix
$\B_n$ of rank $k$ and an error-matrix $\EE_n$ with
$\| \EE_n \| =o (\sqrt{n})$ such that $\A_n =\B_n +\EE_n$ $(n=k,k+1,\dots )$.
Say, $\B_n$ is the blown-up matrix of the $k\times k$ pattern matrix $\PP_n$,
the $ij$  entry $p_{ij}^{(n)}$ of which is the common entry of the 
$U_{in} \times U_{jn}$ block of $\B_n$.


Then using the relation between the cut-norm of a graphon and a matrix,
further, between the cut-norm and
the spectral norm of a matrix, 
and the transformation of a graph into a step-function graphon, we get that
$$
\| W_{\EE_n } \|_{\square} \le \frac1{n^2} \| \EE_n \|_{\square} \le
 \frac1{n^2} n \| \EE_n \| =\frac1{n} o(\sqrt{n}) =o(n^{-1/2}),
$$
where $\| \EE_n \|$ is the spectral-norm, $\| \EE_n \|_{\square}$ is the
matrix cut-norm of $\EE_n$, and  $W_{\EE_n }$ denotes the graphon 
corresponding to the symmetric matrix $\EE_n$ of uniformly bounded entries. 
(The sides of the unit square are 
divided into $n$ equidistant intervals, and the step-function over 
the small squares of the unit
square takes on values corresponding to the  matrix entries.)
Though, these entries can be negative, the theory of bounded graphons applies
to $W_{\EE_n }$, too.

Using the Steiner equality,  
we get that the squared Frobenius norm of \newline $\A_n -\B_n$,
restricted to the $ij$ block, is
$$
\begin{aligned}
 &\| (\A_n -\B_n)_{ij} \|^2_F =\sum_{u\in U_{in}} \sum_{v\in U_{jn} } 
(a_{uv}^{(n)} -p_{ij}^{(n)} )^2 \\
 &=\sum_{u\in U_{in}} \sum_{v\in U_{jn} } 
 (a_{uv}^{(n)} - d(U_{in} ,U_{jn} ) )^2 
 + |U_{in}| |U_{jn}| (d(U_{in} ,U_{jn} )-p_{ij}^{(n)} )^2 ,
\end{aligned}
$$
where the edge-density $d(U_{in} ,U_{jn} )$ of~(\ref{edconv}) 
is now viewed as the average of the entries of $\A_n$
in the $U_{in} \times U_{jn}$ block. 
Then by the inequality between the Frobenius
and spectral norms,
$$
 \| (\A_n -\B_n)_{ij} \|^2_F \le n\| \A_n -\B_n \|^2 = n \|\EE_n \|^2=
  n (o (\sqrt{n}))^2 .
$$
Therefore, for every $1\le i\le j\le k$ pair:
\begin{equation}\label{de}
\begin{aligned}
 (d(U_{in} ,U_{jn} )-p_{ij}^{(n)} )^2 &\le \frac1{|U_{in}| |U_{jn}|}  
n (o (\sqrt{n}))^2 =
\frac1{\frac{|U_{in}|}{n} \frac{|U_{jn}|}{n} }  n (\frac{o (\sqrt{n})}{n})^2  \\
 &=n (o(n^{-1/2}))^2 =o(1)
\end{aligned}
\end{equation}
as $\frac{|U_{in}|}{n} \to r_i$ when $n\to \infty$ $(i=1,\dots ,k)$.
Consequently, $p_{ij}^{(n)} \ge 0$, provided there are no
constantly zero blocks in $\A_n$.

Eventually, we prove the $G_n \to W_H$ convergence by verifying that the
cut-distance between the corresponding graphons tends to 0. 
Using the triangle inequality, we get
$$
\| W_{G_n} -W_{H} \|_{\square} \le
\| W_{G_n} -W_{\B_n } \|_{\square} +
\| W_{\B_n } -W_{G_n /P_{k,n}} \|_{\square}+
\| W_{G_n/P_{k,n}} -W_{H} \|_{\square}
$$
were $G_n /P_{k,n}$ is the factor graph of $G_n$ with respect to the 
$k$-partition $P_{k,n} =(U_{1n} ,\dots , U_{kn})$. 
This is an edge- and vertex-weighted graph on $k$ vertices, with
vertex-weights $\frac{|U_{in}|}{n}$ and edge-weights $d (U_{in}, U_{jn})$,
$i,j=1,\dots ,k$.

The first term is $\| W_{\EE_n } \|_{\square} =o(n^{-1/2})$. 
To estimate the second term, we use that $\B_n$ is the blown-up matrix
of $\PP_n$ with respect to the $k$-partition $P_{k,n}$, after conveniently 
permuting
its rows (and columns, accordingly). The graphon  $W_{\B_n }$ is also stepwise
constant over the unit square, where the sides are divided into $k$ parts:
the interval $I_j$ has lengths $\frac{|U_{jn}|}{n}$ $(j=1,\dots ,k)$, and over 
$I_i \times I_j$ the step-function takes on the value $p_{ij}^{(n)}$.
By its nature, the graphon  $W_{G_n/P_{k,n}}$ is stepwise constant with the
same subdivision of the unit square, and over 
$I_i \times I_j$ it takes on the value $d(U_{in} ,U_{jn} )$, $i,j=1,\dots ,k$.
But in view of (\ref{de}), 
$\| W_{\B_n } -W_{G_n /P_{k,n}} \|_{\square} =\sqrt{n} o(n^{-1/2}) =o(1)$. 
The third term is $o(1)$, because of the assumptions
$\frac{|U_{in}|}{n}\to r_i$ $(i=1,\dots ,k)$ and
$d (U_{in},U_{jn} ) \to p_{ij}$, $i,j=1,\dots ,k$. 
Therefore, $\| W_{G_n} -W_{H} \|_{\square} =o(1)$ and  so,
$G_n \to H$, which finishes the proof. $\square$

\vskip0.5cm

Summarizing: we proved that PIV $\Longleftrightarrow$ P0 $\Longrightarrow$ PI, 
and PI+ $\Longrightarrow$ P0.
So PI is weaker and PI+ is stronger than P0.
We should find something
between PI and PI+ which is equivalent to P0, what is an open question yet. 
Further, \newline P0 $\Longrightarrow$ PII $\Longrightarrow$ PIII. 
We are also able to prove the equivalences 
of Theorem~\ref{ekv} by strengthening PII and PIII.

\newpage
\noindent
\textbf{Proof of Theorem~\ref{ekv}}.
For the proof we need the following lemma.
\begin{lemma}\label{nodom}
Under P0, the following holds for except $o(n_i )$ vertices $u\in U_i$, 
and for every $1\le i \le j \le k$:
$$
 N_1 (u; U_j ) =(1+o(1)) p_{ij} n_j .
$$
\end{lemma}
The statement of the lemma follows from the $P_1 (s)$ $(\forall s ) \Rightarrow
P'_0$  implication of the
Chung--Graham--Wilson~\cite{Chung1} paper (if $i=j$) and its bipartite analogue
(if $i\ne j$).
In the multiclass scenario, Lemma~\ref{nodom} implies 
that there are no dominant vertices, i.e., the vertex-degrees are of
order $\Theta (n)$, except for at most $o(n)$ vertices. Also,
except at most $o(n)$ vertices, the clusters are
asymptotically regular, and the cluster pairs are asymptotically biregular.
Note that this property, even in the $k=1$ case,  is weaker than the equivalent
properties of quasirandomness, in particular, than that for the codegrees. 

\vskip0.2cm
\noindent
\textbf{Proof of P0 $\Longrightarrow$ PII+ $\Longrightarrow$ PIII+:} 
These implications hold in view of
the proof of Theorem~\ref{conj}, as the extra conditions follow from P0,
automatically. Indeed, P0 implies the convergence of the 
edge-densities (see~\cite{LovSos}, the end of the proof of Theorem 2.2) and
that of the within- and between-cluster vertex-degrees (see Lemma~\ref{nodom}),
though, latter is weaker than the quasirandomness. $\square$

\vskip0.2cm
\noindent
\textbf{Proof of PIII+ $\Longrightarrow$ PIV:} 
We will prove that $P_4$,  and the
analogous statement of~\cite{Thomason1} holds, whenever PIII+
holds. Here we cite the notion of $(p-\beta )$-jumbled and bi-jumbled
graphs based on~\cite{Bollobas,Kohayakawa,Thomason,Thomason1}.
A graph is  $(p-\beta )$-jumbled if for any $X\subset V$
$$
 \left| e(X) - p {|X| \choose 2 } \right| \le \beta |X| 
$$
with $0<p\le 1 <\beta$ and $e(X)$ denotes the number of edges of the
underlying graph with both endpoints in $X$. 
In the random graph $G_n (p)$, $\beta = O(\sqrt{n})$, which is best possible
(see~\cite{Bollobas}), and for a quasirandom $G_n$,  $\beta = o(n)$, 
in view of $P_4$ of~\cite{Chung1}.  

Likewise, a bipartite graph is  
$(p-\beta )$-bijumbled if for every $X,Y\subset V$
$$
 \left| e(X,Y) - p |X| |Y | \right| \le \beta \sqrt{|X| |Y|} . 
$$
It makes sense for bipartite graphs with the partition of vertices 
$(U_1 ,U_2 )$, where $X\subset U_1$ and $Y\subset U_2$. For bipartite
quasirandom graphs, $\beta =o(n)$ again, see~\cite{Thomason1}.

Under the degree-conditions, $\Vol (X)$ is approximately $n$ times  $|X|$,
so the right hand side constant, analogous to $\beta$, will be
$o(1)$ in our case, when we use volumes instead of cardinalities.
Indeed, if $i=j$, then we take into consideration that by~(\ref{degri}),
for  $X\subset U_i$,
$$
 \Vol (X) =|X| (1+o(1)) \sum_{\ell =1}^k p_{i\ell} n_ {\ell} .
$$
Let $(U_1 ,\dots ,U_k )$ be the $k$-partition, guaranteed by PIII+,  such that
\newline $\md_k (G_n; U_1 ,\dots U_k ) =o(1)$. Then for $X\subset U_i$, 
$$
\begin{aligned}
 &e(X,X) - p_{ii} |X|^2 =e(X,X) - [d(U_i ,U_i ) +o(1)] |X|^2  \\
& =e(X,X)- \frac{e(U_i ,U_i) }{\frac{\Vol^2 (U_i ) }{(1+o(1))^2 
 ( \sum_{\ell =1}^k p_{i\ell} n_{\ell} )^2 }} \frac{\Vol^2 (X)}
 {(1+o(1))^2 
 ( \sum_{\ell =1}^k p_{i\ell} n_{\ell} )^2 } -o(1) |X|^2  \\
& =[e(X,X) -\rho (U_i ,U_i ) \Vol^2 (X)] -o(1) \rho (U_i ,U_i )\Vol^2 (X)   -
 o(1) |X|^2  \\
&\le \md_k (G_n; U_1 ,\dots U_k ) \sqrt{\Vol^2 (X)} -o(1) e(U_i ,U_i ) 
 \left( \frac{\Vol (X)}{\Vol (U_i )} \right)^2 - o(n^2  ) = o(n^2)
\end{aligned}
$$
as $\md_k (G_n ; U_1 ,\dots U_k ) =o(1)$ by PIII+, 
and we also used~(\ref{densi}). 
Then $P_4$ of~\cite{Chung1} implies $P_7$ of~\cite{Chung1}, that is our PIV.
If $i\ne j$, then with a similar argument, for $X\subset U_i$, $Y\subset U_j$:
$$
\begin{aligned}
 &e(X,Y) - p_{ij} |X| |Y| =e(X,Y) - [d(U_i ,U_j ) +o(1)] |X| |Y| \\
&= e(X,Y) - \frac{e(U_i ,U_j) }{\frac{\Vol (U_i ) \Vol (U_j ) }{(1+o(1))^2 
 ( \sum_{\ell =1}^k p_{i\ell} n_{\ell} ) 
 ( \sum_{\ell =1}^k p_{j\ell} n_{\ell} ) }} \frac{\Vol (X) \Vol (Y) }
 {(1+o(1))^2 
 ( \sum_{\ell =1}^k p_{i\ell} n_{\ell} )
 ( \sum_{\ell =1}^k p_{j\ell} n_{\ell} ) } \\
 &-o(1) |X| |Y|  \\
& =[e(X,Y) -\rho (U_i ,U_i ) \Vol(X) \Vol (Y)] 
 -o(1) \rho (U_i ,U_i )\Vol (X)    \Vol (Y) - o(1) |X| |Y| \\
&\le \md_k (G_n ; U_1 ,\dots U_k ) \sqrt{\Vol(X) \Vol (Y)} -o(1) e(U_i ,U_j ) 
 \left( \frac{\Vol (X)}{\Vol (U_i )}\right) \left( \frac{\Vol (Y)}{\Vol (U_j)}
 \right) \\
 & - o(n^2  ) = o(n^2)
\end{aligned}
$$
as $\md_k (G_n ; U_1 ,\dots U_k ) =o(1)$.
By Theorem 2 of~\cite{Thomason1}, it implies P0, and so, PIV.
Further, by PIII+,  $\md_1 (G_n )>\theta ,\dots , \md_{k-1} (G_n  ) >\theta$,
so $\md_i (G_n ; U'_1 ,\dots ,U'_i ) >\theta$ with any 
$( U'_1 ,\dots ,U'_i )\in {\part}_i$ for every $i=1,\dots ,k-1$.
With the above argument it follows that the requirements of $P_4$ 
of~\cite{Chung1}, see~(\ref{**}), and those of~\cite{Thomason1} 
cannot be met with a 
multiway discrepancy  $\md_i (G_n; U'_1 ,\dots U'_i )$ 
with  $i\in \{ 1,\dots , k-1 \}$ and $(U'_1 ,\dots U'_i ) \in {\part}_i$.
$\square$

\section{Conclusions}\label{concl}

In course of proving equivalences between generalized quasirandom properties,
we have characterized spectra and spectral subspaces of generalized
quasirandom graphs; further, used a version of
the expander mixing lemma to the $k$-cluster case and its certain converse.
These theorems can give a hint for practitioners about the
choice of the number of clusters. The original expander mixing lemma and 
its converse (for simple, regular graphs) treat
the $k=1$ case only, whereas the Szemer\'edi regularity lemma applies to
the worst case scenario: even if there is not an underlying cluster structure,
we can find cluster pairs with small discrepancy with an enormously large $k$
(which does not depend on the number of vertices, it only depends on the
discrepancy to be attained). The lemma also has constructive versions,
and relations to spectra, see Szegedy~\cite{Szegedy}.
Here we rather treat the intermediate case, 
and show
that a moderate $k$ suffices if our graph has $k$ structural eigenvalues 
and a hidden $k$-cluster structure
that can be revealed by spectral clustering tools. Under good clustering
we generally understand clusters with small within- and between-cluster
discrepancies. 


In most of the real-life problems, we want just to find regular partitions 
that minimize 
discrepancies both within and between clusters, such that vertices of the same
cluster behave similarly towards vertices of the same (own or other) cluster.
This new paradigm for structural decomposition (see~\cite{Pelillo})
relies on minimizing the
within- and between-cluster discrepancies by means of the normalized
modularity matrix.
Depending on the
number and sign of the  so-called structural eigenvalues
of this matrix, inferences can be made 
on the number of the underlying clusters and the type
of connection between them.
Furthermore, even if not all of the the generalized quasirandom properties 
hold,
the PII $\Rightarrow$ PIII implication of Theorem~\ref{ekv}
makes it possible to reveal the
structure of the actual data with spectral tools.
The problem can as well be generalized to rectangular arrays, 
e.g., for biclustering genes and 
conditions of microarrays at the same time, when we  want to find clusters of 
similarly functioning genes
that equally (not especially weakly or strongly) influence conditions of the 
same cluster, see~\cite{Bolla10}. 
Sometimes, the graph sequence develops in time (e.g., internet or 
keyword--document matrices), but still we want to find some fixed number
of underlying clusters (like user groups or topics) with spectral clustering.
 
We remark that 
if $r=\rk (\PP ) <k$, then there are only $r$ structural eigenvalues of
the adjacency and normalized modularity matrices of the graph sequence, but the
weighted multiway variance is not necessarily  dropping down at $r$.
If especially, there are identical rows in $\PP$, and
the number of distinct rows is $\ell$ ($r \le \ell \le k$), 
then it is $\md_{\ell}$ and  ${\tilde S}^2_{\ell}$ that drops suddenly 
compared to  $\md_{\ell -1}$ and  ${\tilde S}^2_{\ell -1}$, respectively.
So in practice, when we detect a gap in the spectrum between the $(r-1)$th and
$r$th eigenvalues, it is only an indication of $k>r$ clusters. To find
the optimal $k$, the cluster variances ${\tilde S}^2_{r}$,
${\tilde S}^2_{r+1}$,\dots should be examined until a `small enough'
${\tilde S}^2_{k}$.

\section*{Acknowledgements}

The author thanks Bojan Mohar
for fruitful discussions and collaboration in proving the 
P0 $\Leftrightarrow$ PIV
equivalence. We are also indebted to Vera T. S\'os and L\'aszl\'o Lov\'asz
for useful suggestions. 
 The research reported in this paper was supported by the Higher 
Education Excellence Program of the Ministry of Human Capacities in the 
frame of  Artificial Intelligence research area of Budapest University 
of Technology (BME FIKP-MI/SC).
We also thank  the VTT Technical Research Centre of Finland for considering
the involved methods in their STOMOGRAPH project.


\begin{thebibliography}{99}

\bibitem{Abbe}
Abbe, E., Sandon, C., Community detection in general stochastic block models:
fundamental limits and efficient algorithms for recovery,
2015 IEEE 56th Annual Symposium on FOCS.

\bibitem{Alon0}
Alon, N., 1986 Eigenvalues and expanders, \textit{Combinatorica}
\textbf{6} (1986), 83-96.


\bibitem{Alon1}
Alon, N., Krivelevich, M. and Vu, V. H., On the concentration of 
eigenvalues of random symmetric matrices, \textit{Isr. J. Math.}
\textbf{131} (2002), 259--267.

\bibitem{Alon10}
Alon, N., Coja-Oghlan, A., Han, H., Kang, M., R\"odl, V. and Schacht, M.,
Quasi-randomness and algorithmic regularity for graphs with general
degree distributions, \textit{Siam J. Comput.} \textbf{39} (2010), 2336--2362.

\bibitem{Azran} 
Azran, A. and Ghahramani, Z., Spectral methods for automatic
multiscale data clustering. In \textit{Proc. IEEE Computer
Society Conference on Computer Vision and Pattern Recognition (CVPR 2006),
New York NY} (Fitzgibbon A, Taylor CJ and Lecun Y eds), 
IEEE Computer Society, Los Alamitos, California (2006),
pp.~190--197. 

\bibitem{Bolla5}
Bolla, M., Recognizing linear structure in noisy matrices,
\textit{Linear Algebra Appl.} \textbf{402} (2005), 228-244.

\bibitem{Bolla10}
Bolla, M., Friedl, K. and Kr\'amli, A., Singular value decomposition
of large random matrices (for two-way classification of microarrays),
\textit{J. Multivariate Anal.} \textbf{101} (2010), 434--446.

\bibitem{Bolla11}
Bolla, M., Penalized versions of the Newman--Girvan modularity and their
relation to normalized cuts and k-means clustering,
\textit{Phys. Rev. E} \textbf{84}  (2011),
016108. 

\bibitem{Bolla12} 
Bolla, M., K\'oi, T. and Kr\'amli, A., Testability of
minimum balanced multiway cut densities, \textit{Discret. Appl. Math.}
\textbf{160} (2012), 1019--1027.


\bibitem{Bolla13}
Bolla, M., Spectral clustering and biclustering. Wiley (2013).

\bibitem{Bolla14}
Bolla, M., Modularity spectra, eigen-subspaces and structure of weighted
graphs, \textit{European J. Combin.} \textbf{35} (2014), 
105--116. 


\bibitem{Bolla16}
Bolla, M., Relating multiway discrepancy and singular values of 
nonnegative rectangular matrices, 
\textit{Discrete Applied Mathematics} \textbf{203} (2016), 26--34.

\bibitem{Bolla17}
Bolla, M., Elbanna, A., Discrepancy minimizing spectral clustering, 
\textit{Discrete Applied Mathematics} \textbf{243} (2018), 286--289.








\bibitem{Bollobas}
Bollob\'as, B., \textit{Random Graphs}, 2nd edn. Cambridge Univ.
Press, Cambridge (2001).


\bibitem{Bollobas3}
Bollob\'as, B., Janson, S. and Riordan, O., The phase transition in
inhomogeneous random graphs, \textit{Random Struct. Algorithms}
\textbf{31} (2007), 3--122.

\bibitem{Borgs}
Borgs, C, Chayes, J. T., Lov\'asz, L., T.-S\'os, V. and
Vesztergombi, K., Convergent graph sequences I: Subgraph
Frequencies, metric properties, and testing, \textit{Advances in Math.}
\textbf{219} (2008), 1801--1851.

\bibitem{Borgs3} 
Borgs, C, Chayes, J. T., Lov\'asz, L., T.-S\'os, V. and 
Vesztergombi, K., Convergent sequences of dense graphs II: Multiway
cuts and statistical physics,
\textit{Ann. Math.} \textbf{176} (2012), 151--219.



\bibitem{Cheeger}
Cheeger, J.,  A lower bound for the smallest eigenvalue of the Laplacian,
In \textit{Problems in Analysis} (ed. R. C. Gunning),
Princeton Univ. Press, Princeton NJ (1970),  pp.~195--199.
paper dedicated to Salmon Bochner

\bibitem{Chung1}
Chung, F., Graham, R. and Wilson, R. K., Quasi-random graphs,
\textit{Combinatorica} \textbf{9} (1989), 345--362.

\bibitem{Chung2}
Chung, F. and Graham, R., Quasi-random graphs with given degree
sequences, \textit{Random Struct. Algorithms} \textbf{12} (2008), 1--19.






\bibitem{Cojab}
Coja-Oghlan, A. and Lanka, A., Finding planted partitions in random graphs
with general  degree distributions. \textit{J. Discret. Math.}
\textbf{23} (2009), 1682--1714.




\bibitem{Diaconis}
Diaconis, P. and Stroock, D., Geometric bounds for eigenvalues of Markov
chains. \textit{Ann. Appl. Probab.} \textbf{1} (1991), 36--62.

\bibitem{Evra}
Evra, S., Golubev, K., Lubotzky, A., Mixing properties and the chromatic
number of Ramanujan complexes, arXiv:1407.7700 [math.CO] (2014).

\bibitem{Fiedler}
Fiedler, M., 1973 Algebraic connectivity of graphs, \textit{Czech. Math. J.}
\textbf{23} (1973), 298--305.


\bibitem{Frieze}
Frieze, A. and Kannan, R., Quick approximation to matrices
and applications, \textit{Combinatorica} \textbf{19} (1999), 
175--220.


\bibitem{Hoffman2}
Hoffman, A. J., Eigenvalues and partitionings of the edges of a graph,
\textit{Linear Algebra Appl.} \textbf{5} (1972), 137--146.

\bibitem{Holland}
Holland, P., Laskey, K. B. and Leinhardt, S., Stochastic blockmodels: 
some first steps. \textit{Social Networks} \textbf{5} (1983), 109--137.


\bibitem{Hoory}
Hoory, S., Linial, N. and Widgerson, A., Expander graphs
and their applications, \textit{Bull. Amer. Math. Soc. (N. S.)}
\textbf{43} (2006), 439--561.



\bibitem{Kohayakawa}
Kohayakawa, Y., R\"odl, V., Schacht, M., Sissokho, P., Skokan, J., 
Tur\'an's theorem for pseudo-random graphs,
\textit{J. Comb. Theory, Ser. A} \textbf{114} (2007), 631--657.

\bibitem{Lee}
Lee, J. R., Gharan, S. O. and Trevisan, L., Multi-way spectral partitioning and
higher-order Cheeger inequalities. In \textit{Proc. 44th Annual 
ACM Symposium on the Theory of Computing (STOC 2012)}, pp.~1117--1130. 
New York NY (2012).


\bibitem{Lovasz0} 
Lov\'{a}sz, L., Random walks on graphs: a survey.
In \textit{Combinatorics, Paul Erd\H os is Eighty. J\'anos Bolyai Society, 
Mathematical Studies} Vol. 2, pp.~1--46. Keszthely, Hungary (1993). 


\bibitem{LovSos}
Lov\'asz, L. and T.-S\'os V., Generalized quasirandom
graphs, \textit{J. Comb. Theory B} \textbf{98} (2008), 146--163.


\bibitem{McSherry}
McSherry, F., Spectral partitioning of random graphs. In  \textit{Proc.
42nd Annual Symposium on Foundations of Computer Science
(FOCS 2001)}, pp.~529--537. Las Vegas, Nevada (2001).

\bibitem{Meila}  
Meila, M. and Shi, J., Learning segmentation by random walks.
In \textit{Proc. 13th Neural Information
Processing Systems Conference (NIPS 2001)} (Leen TK,
Dietterich TG and Tresp V eds),  pp.~873--879. MIT Press, Cambridge, USA
(2001).

\bibitem{Mohar1}
Mohar, B., Isoperimetric inequalities, growth and the spectrum of graphs. 
\textit{Linear Algebra Appl.} \textbf{103} (1988), 119--131.

\bibitem{Mossel}
Mossel, E., Neeman, J. and Sly, A., Reconstruction and estimation in the
planted partition model, \textit{Probab. Theory Related Fields} \textbf{162}
(2015), 431--461.

\bibitem{Newman}
Newman, M. E. J., \textit{Networks, An Introduction}. Oxford University Press
(2010).

\bibitem{Ng} 
Ng, A. Y., Jordan, M. I. and Weiss, Y.,
On spectral clustering: analysis and an algorithm.
In \textit{Proc. 14th Neural Information
Processing Systems Conference (NIPS 2001)} (Dietterich TG,
Becker S and Ghahramani Z eds), pp.~849--856. MIT Press, Cambridge, USA (2001).
Vancouver

\bibitem{Ostrovsky}
Ostrovsky et al., The effectiveness of Lloyd-type methods for the $k$-means
problem, \textit{J. ACM} \textbf{59} (6), Article 28 (2012).

\bibitem{Pelillo}
Pelillo, M., Elezi, I., Fiorucci, M., Revealing Structure in
Large Graphs: Szemer\'edi's Regularity Lemma and Its Use in Pattern
Recognition, \textit{Patter Recognition Letters} \textbf{87} (2017),
4--11. 


\bibitem{SimonovitsS}
Simonovits, M. and  T.-S\'os, V., Szemer\'edi's partition and
quasi-randomness, \textit{Random Struct. Algorithms} \textbf{2} (1991), 1--10.

\bibitem{Szegedy}
Szegedy, B., Limits of kernel operators and the spectral regularity lemma,
\textit{European J. Combin.} \textbf{32} (2011), 
1156-1167.

\bibitem{Szemeredi}
Szemer\'edi, E., Regular partitions of graphs. In \textit{Colloque
Inter. CNRS. No. 260, Probl\'emes Combinatoires et
Th\'eorie Graphes} (Bermond J-C, Fournier J-C, Las Vergnas M and
Sotteau D eds), pp.~399--401 (1976).

\bibitem{Thomason}
Thomason, A., Pseudo-random graphs, \textit{Ann. Discret. Math.} \textbf{33}
(1987), 307--331.

\bibitem{Thomason1}
Thomason, A., Dense expanders and pseudo-random bipartite graphs,
\textit{Discret. Math.} \textbf{75} (1989), 381--386.


\end{thebibliography}
\end{document}